\documentclass{amsart} 
\usepackage{amsfonts,amsmath,amsthm,amssymb}
\newtheorem{lemma}{Lemma}
\newtheorem{theorem}{Theorem}
\newtheorem{conjecture}{Conjecture}
 
\title{Symmetric subgroups in modular group algebras}

\author{A.B.~Konovalov, A.G.~Krivokhata}

\date{}

\address{
A.B.~Konovalov
\newline School of Computer Science, University of St Andrews,
\newline Jack Cole Building, North Haugh, St Andrews, Fife, KY16 9SX, Scotland}
\email{konovalov@member.ams.org}

\address{
A.G.~Krivokhata (nee Tsapok)
\newline Department of Mathematics, Zaporozhye National University, 
\newline 66 Zhukovskogo str., 69063, Zaporozhye, Ukraine}
\email{k-algebra@zsu.zp.ua}



\begin{document}

\begin{abstract}
Let $V(KG)$ be a normalised unit group of the modular group algebra
of a finite $p$-group $G$ over the field $K$ of $p$ elements.
We introduce a notion of symmetric subgroups in $V(KG)$ as subgroups
invariant under the action of the classical involution of the group algebra
$KG$. We study properties of symmetric subgroups and construct a
counterexample to the conjecture by V.~Bovdi, which states that
$V(KG)=\langle G, S_* \rangle$, where $S_*$ is a set of symmetric
units of $V(KG)$.
\end{abstract}

\maketitle

\section{Introduction}

\noindent In this preprint\footnote{This preprint is a slightly 
updated translation of the paper \cite{Kon-Tsap}, earlier
written by the same authors in Russian.} we investigate the 
structure of the unit group $U(KG)$ of the modular group algebra
$KG$, where $G$ is a finite $p$-group and $K$ is a field of
$p$ elements. $U(KG)$ is a direct product of the unit group
of the field $K$ and the normalized unit group $V(KG)$ that
consists of all elements of the form $1+x$, where $x$ belongs
to $I(G)$, the latter being the augmentation ideal of the 
modular group algebra $KG$. Thus, the task of the investigation 
of $U(KG)$ is reduced to the investigation of the $p$-group 
$V(KG)$.

Study of units and their properties is one of the main
research problems in group ring theory. Results obtained
in this direction are also useful for the investigation of Lie
properties of group rings, isomorphism problem and other
open questions in this area (see, for example, \cite{Bovdi-Kurdics}).

Let us consider the mapping $\varphi$ from $KG$ into itselfâ 
which maps an element 
$x = \sum \lambda_g g$ \; to \; $\sum \lambda_{g} g^{-1}$. 
This mapping is an antiautomorphism
of order two, and it is called the {\it classical involution} of
the group algebra $KG$. In what follows, we will denote $\varphi(x)$
by $x^*$. An element $x \in KG$ will be called {\it symmetric}, if
$x=x^*$. Note that if $x^* \neq x^{-1}$, then $x^*x$ --- nontrivial
symmetric element. Indeed, $\varphi(x^*x)=(x^*x)^* = x^* x^{**} = x^*x$
since $x^{**}=x$. Classical involution and symmetric elements 
were studied by V.~Bovdi, M.~Dokuchaev, L.~Kovacs, S.~Sehgal in 
\cite{BKS, Bovdi-SU, Bovdi-Dok, Bovdi-Roz, Sehgal-Sym}.

In \cite{BovdiAA-units} the generating system for $V(KG)$ was given,
however that system is minimal only in the abelian case. In the general case,
the problem of the determination of the minimal generating system of
$V(KG)$ remains open. For groups of small order the minimal generating
system can be determined with the help of computer after the construction
of $V(KG)$ using the algorithm from \cite{BovdiAA-units} implemented in
the LAGUNA package \cite {LAGUNA} for the computational algebra system
GAP \cite {GAP}. Obtaining new results about the minimal generating 
system would be useful for the improvement of algorithms for the
computation of $V(KG)$.

In 1996 V.~Bovdi suggested the following conjecture:

\begin{conjecture}
Let $G$ be a finite nonabelian $p$-group, $V(KG)$ be the normalized
unit group of the modular group algebra $KG$, and 
$S_* = \{ x \in V(KG) \ | \ x^* = x \}$ be the set of symmetric units.
Is it true that
\begin{equation}
\label{gip-Bovdi} V(KG)=\langle G, S_* \rangle ?
\end{equation}
\end{conjecture}

In the present paper we give a counterexample to this conjecture
and investigate further properties of the classical involution.
We introduce the notions of symmetric subsets of the group algebra
and symmetric subgroups of $V(KG)$, and investigate their properties.
Then we give some obvious examples of symmetric subgroups and formulate
the question about the existence of non-trivial symmetric subgroup.
After this, we consider the normalized unit group of the modular group
algebra of the quaternion group of order 8, and check for it the 
conjecture (\ref{gip-Bovdi}). Using the LAGUNA package \cite{LAGUNA}, we
discovered that in this case the condition (\ref{gip-Bovdi}) does not
hold, because $\langle G, S_* \rangle$ has the order 64. In the present 
paper we give a purely theoretical proof of this fact. Besides this,
it appears that $\langle G, S_* \rangle$ is the symmetric subgroup,
as well as the set $S_*$, which in this case is also a subgroup. Thus,
we obtain two examples of non-trivial symmetric subgroups.

It would be interesting to continue studies of properties of symmetric
subgroups of $V(KG)$ and conditions of the existence of non-trivial 
symmetric subgroups in $V(KG)$.

\section{Symmetric subgroups}

\noindent Let $H$ be a subset of the group algebra $KG$. 
Then $H$ will be called {\it symmetric}, if $H^*=H$, where 
$H^* =\{ h^* | h \in H \}$. Similarly, a subgroup 
$H \subseteq V(KG)$ will be called a {\it symmetric subgroup},
if $H^*=H$.

Clearly, symmetric subgroups exist. It is easy to see that
$\{1\}$, $G$ and $V(KG)$ are trivial examples of symmetric
subgroups. This naturally raises a question of the existence
of modular group algebras that possesses non-trivial symmetric
subgroups. Moreover, in \cite{BKS} conditions were obtained
under which the set of symmetric elements forms a group,
that in this case will be symmetric. Therefore, it would be 
interesting to find an example of non-trivial symmetric subgroup
that contains non-symmetric units, thus the restriction of
the classical involution on this subgroup will be not an identity 
mapping.

Before we construct such example, we will state and prove some
properties of symmetric subgroup.

\begin{lemma}
\label{lemma1}
If $H$ is a subgroup of $V(KG)$, then
$H^*$ also is a subgroup of $V(KG)$.
\end{lemma}

{\it Proof.} We will check that $H^*$ is a subgroup of $V(KG)$.
Since $1 \in H$, we have $1=\varphi(1) \in H^*$. Then, if 
$a^*,b^* \in H^*$, so does $a^* b^* \in H^*$, since 
$a^*b^*=\varphi(ba)$, and $ba \in H$. Finally, if $a^* \in H^*$, 
then $(a^*)^{-1}=(\varphi(a))^{-1}=\varphi(a^{-1}) \in H^*$, and
the lemma is proved.

\begin{lemma}
\label{lemma2}
Let $H'$ and $H$ be subgroups of $V(KG)$, and $H$
be a symmetric subgroup. If $H'$ is a conjugate of $H$ with the
help of an element $g \in V(KG)$ such that $(g^*)^{-1}=(g^{-1})^*$, 
then $H'^*$ is also a conjugate of $H$ with the help of element 
$(g^*)^{-1}$.
\end{lemma}

{\it Proof.} Let us assume that $g \in V(KG)$ is such that
$H'=g^{-1}Hg$, and $H=H^*$. Then 
$H'^*=(g^{-1}hg)^*=g^*h(g^{-1})^*=g^*H(g^{-1})^*$. 
If $(g^*)^{-1}=(g^{-1})^*$, what holds, in particular, for all
$g \in G$, then $H'^*=H^{(g^{-1})^*}$, and the lemma is proved.
Note that $H'$ is not necessary to be symmetric, and also that
$H'^*=(H^{g^{-1}})^*$.

\begin{lemma}
\label{lemma3}
Let $H_1$ and $H_2$ be symmetric subgroups.
Then $H_1 H_2$ also is a symmetric subgroup.
\end{lemma}

{\it Proof.} Let $H_1$ and $H_2$ be symmetric subgroups.
Consider an element from $H_1 H_2$ having the form 
$h = a_1 b_1 a_2 b_2 \dots a_k b_k$, where 
$a_i \in H_1$ and $b_i \in H_2$. Then 
$h^* = u_k v_k \dots u_2 v_2 u_1 v_1$, where $u_i=b_i^*, v_i=a_i^*$, 
also belongs to $H_1 H_2$, because $H_1$ and $H_2$ are symmetric.
Thus, the lemma is proved.

\section{Investigation of the conjecture on generators of $V(KG)$}

\noindent To verify the conjecture by V.~Bovdi about generators of
the normalized unit group we used the LAGUNA package \cite{LAGUNA} 
for the computational algebra system GAP \cite{GAP}. We developed a 
program in GAP language using functions available in the LAGUNA
package, that performed the following steps:

\begin{list}{\alph{enumi})}{}\usecounter{enumi}

\item
get the set of all elements of the group $G$;

\item
compute the normalized unit group $V(KG)$ 
(generated by group algebra elements) with the help of 
the LAGUNA package, and the group $W$ that is a 
power-commutator presentation of $V(KG)$;

\item
construct the embedding $f$ from $G$ 
to the group algebra $KG$;

\item
construct the isomorphism $\psi$ 
from $V(KG)$ to the group $W$;

\item
compute the set $S$ of symmetric units of the group
algebra $KG$;

\item
using the isomorphism $\psi$, find images of
the set of symmetric units $S$ and $f(G)$ in the
group $W$;

\item
compute the subgroup $\langle \psi(f(G)), \psi(S_*) \rangle \subseteq W$
and verify the condition (\ref{gip-Bovdi}).

\end{list}

As a result, we
find out that the condition (\ref{gip-Bovdi}) does not hold for the
quaternion group of order eight.

\section{Computation of the symmetric subgroup $\langle G,S_* \rangle$ }

\noindent In this section we will give a theoretical explanation
why the quaternion group $Q_8$ provides a counterexample to
the conjecture by V.~Bovdi. The group $Q_8$ is given by the
following presentation:
$$
Q_8=\langle a,b | a^4=b^4=1, a^2=b^2, b^{-1}ab=a^3 \rangle,
$$
and consists of the following elements:
$$
Q_8=\{a^ib^j, 0\leq i \leq 4, 0\leq j\leq 2\}.
$$
Every symmetric element $x \in S_*$ has the form
\begin{equation}
\label{eks} x=\alpha_0+\alpha_1a^2+\gamma_x,
\end{equation}
where
\begin{equation}
\label{gamma}\gamma_x=\alpha_2(a+a^3)+\alpha_3(b+a^2b)+\alpha_4(ab+a^3b).
\end{equation}
Since $x$ is a unit from $V(KQ_8)$, we have that
$\alpha_0+\alpha_1=1$ (i.e. one and only one of coefficients
$\alpha_0, \alpha_1$ is equal to one). From this follows that
$|S_*|=16$.

\begin{lemma}
\label{lemma4}
Let $S_*=\{s \in V(KQ_8) | s^*=s\}$ be the set of
symmetric units. Then $S_*\subset Z(KQ_8)$, where $Z(KQ_8)$ is the
center of the group algebra $KQ_8$.
\end{lemma}

{\it Proof.} The group $Q_8$ has the following decomposition
on the conjugacy classes of elements:
$$
Q_8 = \{ 1 \} \cup \{ a^2 \} \cup \{ a, a^3 \} \cup 
      \{ b, a^2b \} \cup \{ ab, a^3 \}.
$$
Since class sums are central elements (see \cite{BovdiAA-GR}),
the lemma follows immediately from (\ref{eks}) and (\ref{gamma}).

\begin{lemma}
\label{lemma5} 
The set $S_*$ is a subgroup of the normalized
unit group $V(KQ_8)$.
\end{lemma}

{\it Proof.} The lemma follows from the result from \cite {BKS} 
stating that the set of symmetric units is a subgroup if and only 
if they all commute, in combination with the lemma \ref{lemma4}.

It also follows from the lemma \ref{lemma5} that $S_*$ is a non-trivial
symmetric subgroup of $V(KQ_8)$. However, all its elements are
symmetric, so we did not have an example of a symmetric
subgroup which contain non-symmetric units yet.

\begin{lemma}
\label{lemma6} 
$\forall x\in S_*$ \quad $x^2=1$.
\end{lemma}

{\it Proof.} By (\ref{eks}), we have that
$$
x^2=\alpha_0^2+\alpha_1a^4+\gamma_x^2=\alpha_0^2+\alpha_1^2=1,
$$
since the characteristic of the field $K$ is two,
elements $1, a^4$ and $\gamma_x$ are in the center 
$Z(KQ_8)$, and $\gamma_x^2=0$.

Now we are ready to prove the main statement, that gives
a counterexample to the conjecture by V.~Bovdi.

\begin{theorem}
Let $Q_8$ be the quaternion group of order eight,
$V(KQ_8)$ be the normalized unit group of its modular group algebra
over the field of two elements, and $S_*=\{x \in V(KQ_8) | x^*=x\}$ be
the set of symmetric units. Then
$$
|H|=| \langle Q_8,S_* \rangle | = 64. 
$$
\end{theorem}

{\it Proof.} Since $S_*\subset Z(KQ_8)$, then each element from
$H$ can be written as $x=gs$, where $g \in Q_8$ and $s \in S_*$.
Note that $Q_8 \bigcap S_*=\{1, a^2\}$ that coincides with $Z(Q_8)$. 
Thus, $\langle Q_8, S_* \rangle$ is a central product of groups
$Q_8$ and $S_*$. From this follows that 
$|H|={{|S_*|\cdot |Q_8|} \over {|Z(Q_8)|}}={{16\cdot 8} \over
{2}}=64$, and the theorem is proved. 

Thus, we proved that the order of a subgroup generated by the group
$Q_8$ and the set of symmetric elements $S_*$ is equal to $64$, so it
does not coincide with the normalized unit group $V(KQ_8)$ of order
$128$. So, $H \neq V(KQ_8)$, and the quaternion group of order eight
gives a counterexample to the conjecture by V.~Bovdi. Besides this,
$H = \langle Q_8, S_* \rangle$ is a non-trivial symmetric subgroup
of the group $V(KQ_8)$, containing both symmetric and non-symmetric
units.


\end{document}